\documentclass{elsarticle}
\usepackage[colorlinks]{hyperref}
\usepackage{amsmath}
\usepackage{amsfonts}
\usepackage{amsthm}
\usepackage{amssymb}
\usepackage{graphicx}
\usepackage{float}
\usepackage{caption}
\usepackage{wrapfig}

\usepackage{amssymb,amsmath,graphicx,psfrag,amsfonts,latexsym,amsthm,amscd}

\newtheorem{theorem}{Theorem}

\newtheorem{prop}[theorem]{Proposition}

\theoremstyle{definition}

\let\today\relax
\makeatletter
\def\ps@pprintTitle{%
    \let\@oddhead\@empty
    \let\@evenhead\@empty
    \def\@oddfoot{\footnotesize\itshape
         {Submitted preprint} \hfill\today}%
    \let\@evenfoot\@oddfoot
    }
\makeatother

\begin{document}

\title{Conflict graphs of maximally planar subgraphs of Petersen family graphs}


\author[a1]{Joel Foisy\corref{c1}}

\cortext[c1]{Corresponding author}
\ead{foisyjs@potsdam.edu}

\author[a1]{Justin Raimondi}
\ead{justinraimondi@gmail.com}
\address[a1]{Department of Mathematics, SUNY Potsdam, Potsdam, NY 13676}


\begin{abstract}
The purpose of this paper is to show that all maximally planar subgraphs of graphs in the Petersen Family have associated conflict graphs unbalanced. All but three strong conflict graphs arising from Petersen Family Graphs are unbalanced, and the three that are balanced all come from $K_{4,4}-e$.  \\
\end{abstract}
\maketitle


This paper consists of details that complement the paper \cite{foisy}. For more context and definitions, please see that paper.\\

We will use the following well-known facts about planar graphs:  
\begin{enumerate}

\item A maximal planar subgraph of a connected graph is connected.
\item Given an $S^2$ embedding of a 2-connected graph, every region is bounded by a cycle (see, for example, \cite{west}).
\item Let $G$ be a graph and $M$ a subgraph of  $G$ that is maximally planar ($M$ is planar, but the insertion of any edge from $G-M$ into $M$ results in a nonplanar graph). Then $M$ is a spanning subgraph of $G$. Equivalently, the only possible $M-$fragments are individual edges.
\end{enumerate}

\begin{figure}

\vskip -.5in
\hskip .8in \includegraphics[scale=.33]{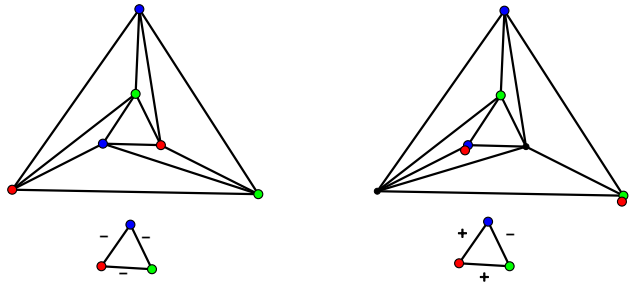}
\vskip -.65in
\caption{The two maximal planar subgraphs of $K_6$.}
\label{fig:k6}
\end{figure}

\begin{prop}
There are two maximal planar subgraphs of $K_6$ and the strong conflict graph for both is not balanced. These subgraphs are pictured in Figure \ref{fig:k6}
\end{prop}

\begin{proof} The first maximal planar subgraph pictured results from removing three non-adjacent edges, the second results from removing a path of length 3. We argue that these are the only maximal planar subgraphs of $K_6$ (up to symmetry). Consider the edges that are removed from $K_6$ to form a maximal planar subgraph. Denote by $H$ the subgraph formed by these edges. If $H$ is connected, then $H$ must contain a path of length 3, else all edges of $H$ are incident to a common vertex, which implies that $K_6-H$ contains $K_5$, and thus is non-planar. Thus, if $H$ is connected, the resulting maximal planar subgraph is the second one pictured in Figure \ref{fig:k6}.

 Suppose $H$ is not connected. If $H$ has 3 components, then the complement of $H$ is the first maximal planar subgraph in Figure \ref{fig:k6}. Suppose $H$ has two components. If each component contains 3 vertices, then $K_6-H$ contains $K_{3,3}$ as a subgraph and is non-planar. Otherwise, one component contains 4 vertices and the other contains 2. The component with 4 vertices cannot contain 2 disjoint edges, else $K_6-H$ is a subgraph of the first subgraph pictured in Figure \ref{fig:k6}. It also cannot contain a triangle, as this would imply that it also contained 2 disjoint edges. It follows that the component with 4 vertices must be a $Y$, but then $K_6-H$ contains a $K_5$ minor. Thus, there is no $H$ with exactly 2 components.
\end{proof}

\begin{figure}
\hskip .15in\includegraphics[scale=.5]{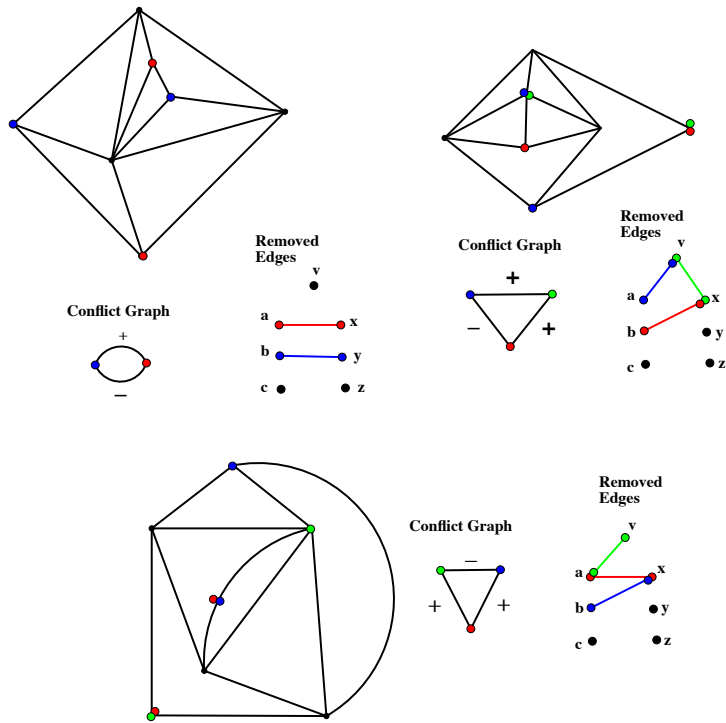}
\vskip -.85in
\caption{The three maximal planar subgraphs of $K_{3,3,1}$.}
\label{fig:k331}
\end{figure}

\begin{prop}
There are three maximal planar subgraphs of $K_{3,3,1}$, and the strong conflict graph for each is not balanced. These subgraphs are pictured in Figure \ref{fig:k331}.
\label{prop:k331}
\end{prop}
\begin{proof}
We partition the vertices of $K_{3,3,1}$ into $\{v\}$, $\{a,b,c\}$ and $\{x,y,z\}$. The first maximal planar subgraph in Figure \ref{fig:k331} results from removing two non-adjacent edges between $\{a,b,c\}$ and $\{x,y,z\}$. No nonisomorphic maximal planar subgraph results from the removal of two such edges. Suppose another maximally planar subgraph results from removing exactly two edges in $G-v$, both of which are incident to the same vertex. Without loss of generality, suppose that an edge is removed from $a$ to $x$ and from $b$ to $x$. Then we must also remove an edge between $v$ and $\{a,b,c\}$, else the non-planar subgraph induced by $\{v,y,z\}$ and $\{a,b,c\}$ is present. It follows that either the third graph pictured in Figure \ref{fig:k331} results (up to symmetry), or that the edges $(v,a)$ and $(v,b)$ are present, but the edge $(v,c)$ is removed. If edge $(v,x)$ is removed, we get a graph equivalent to a subgraph of the second graph in Figure \ref{fig:k331}.  Otherwise, as the edge $(v,a)$ is present, then by contracting, we get a graph that contracts onto $K_{3,3}$, by contracting the edge $(v,a)$, and making one partition $\{(v,a),b,c\}$ and the other $\{x,y,z\}$. 

If more than three edges are removed from $K_{3,3,1}-v$, then we obtain a graph equivalent to a subgraph of the first graph in Figure \ref{fig:k331}. Now suppose exactly three edges have been removed from $K_{3,3,1}-v$, with all incident to a common vertex (say, $x$). In order for the resulting graph to be planar, an edge from $v$ to $\{a,b,c\}$ must also be removed. This results in a graph isomorphic to a subgraph of the third graph pictured in Figure \ref{fig:k331}.

Finally, suppose that just one edge has been removed from $K_{3,3,1}-v$. Without loss of generality, we may assume this edge is $(b,x)$. If $v$ is connected to $x$ and to $b$, then the resulting graph is non-planar, even if other edges are removed from $v$. Without loss of generality, suppose the edge $(v,x)$ has been removed. If none of the edges $(v,a),(v,b)$ or $(v,c)$ is removed, the resulting graph contains a $K_{3,3}$, by partitioning $\{a,b,c\}$ and $\{x,y,z\}$. If $(v,a)$ or $(v,c)$ is removed, we obtain a subgraph of the second graph in Figure \ref{fig:k331}. Otherwise, if only the additional edge $(v,b)$ is removed, then the result is $K_{3,3,1}-\Delta$, which is non-planar. Moreover, any additional removed edge will result in a subgraph of the second graph pictured in Figure \ref{fig:k331}. This concludes our proof. 

One can check (see Figure \ref{fig:k331}) that the conflict graphs of each maximally planar subgraph is balanced.
\end{proof}

\vskip -.75in

\vskip .75in 
Before we discuss the maximal planar subgraphs of $G_8$, we introduce the following:

\begin{prop} Suppose $H$ is obtained from $G$ by a $\Delta-Y$ exchange. Let $S$ represent the set of maximal planar subgraphs of $G$ that contain the triangle used in the $\Delta-Y$ exchange. Let $S'$ represent the set of graphs obtained by a $\Delta-Y$ operation on an element from $S$. Then any maximally planar subgraph of $H$ that contains the $Y$ must be a subgraph of an element of $S'$.
\label{prop:delta}
\end{prop}

\begin{proof}
The result follows from the well-known fact that the $Y-\Delta$ operation preserves planarity. Let $P'$ be a maximally planar subgraph of $H$ that contains the $Y$. By applying a $Y-\Delta$ operation on $P'$, we obtain a planar subgraph (say $P$) of $G$ that contains the triangle. Furthermore, $P$ must be a subgraph of an element of $S$, as $S$ contains every maximally planar subgraph of $G$. It follows that $P'$ must be a subgraph of an element of $S'$.

\end{proof}

\begin{prop} The fourteen maximally planar subgraphs of $G_8$ represented in Figure \ref{fig:p8} and Figure \ref{fig:p8b} represent all possible maximally planar subgraphs of $G_8$. Each has an unbalanced strong conflict graph.
\end{prop}

\begin{proof}
We first note that in Figure \ref{fig:p8}, the vertex $w$ represents the $Y$ vertex  of $G_8$ that is obtained from doing a $\Delta-Y$ exchange on $K_{3,3,1}$. The vertex $v$ represents the vertex of degree $5$ in $G_8$, and $\{a,b,c\}$ and $\{x,y,z\}$ represent the two partitions of size three in $K_{3,3,1}$ that would result if the $Y$ is swapped back to a triangle. The colored edges in Figure \ref{fig:p8} represent edges of $G_8$ that are removed to obtain the maximal planar subgraph. If a maximal planar subgraph of $G_8$ contains all of the $Y$ edges, then by Proposition \ref{prop:delta} and Proposition \ref{prop:k331} it must be a subgraph of one of the second through seventh graphs pictured in Figure \ref{fig:p8}, or the first graph, with edge $(a,y)$ put back in. We denote the first graph with edge $(a,y)$ put back in $G'$. It can be verified (see Figure \ref{fig:p8second}) that the second through eighth graphs are planar, and thus they must be maximally planar. The graph $G'$ requires the removal of an additional edge to obtain a planar graph. If either edge $(b,x)$ or $(a,y)$ is removed, we get the same (isomorphic) maximally planar subgraph of $G_8$, the first pictured in Figure \ref{fig:p8}. If edge $(b,z)$ or $(c,y)$ or $(c,x)$ or $(z,a)$ is removed, we obtain a subgraph of the second graph pictured. If an edge with $v$ as an endpoint is removed, note that if both $(v,b)$ and $(v,y)$ are removed, we still have a $K_{3,3}$ minor on $\{v,c,z\}$ and $\{w,(b,x), (a,y)\}$. As there are no other non-$Y$ edges to remove, we have listed all possible maximal planar subgraphs of $G_8$ obtained from subgraphs of $G'$.
\vskip -1.3in
\begin{figure}
\hskip .3 in\includegraphics[scale=.43]{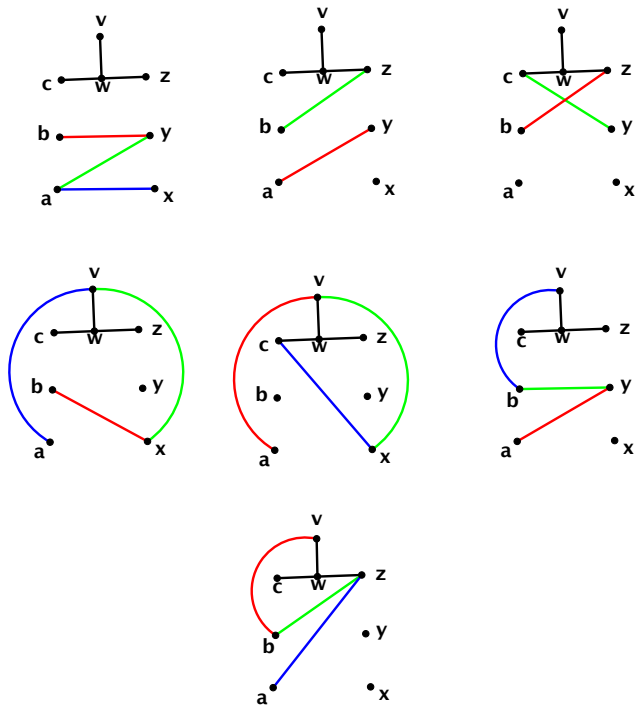}
\vskip -.5in
\caption{The representations of the seven maximal planar subgraphs of $G_8$ that include all edges of the $Y$ that would transform it back to $K_{3,3,1}$.}
\label{fig:p8}
\end{figure}

\vskip 1.3in
Now we consider maximal planar subgraphs that result from removing one or more edges of the $Y$. We first consider exactly one edge being removed from the $Y$. Suppose the edge removed is incident to $v$ (that is, edge $(v,w)$). A maximal planar subgraph that has only edge $(v,w)$ removed from the $Y$ must have a non $Y$ edge removed. First suppose no edge equivalent to $(b,y)$ has been removed. An edge not incident to $v$ must be removed, so we may assume, without loss of generality, that edge $(b,z)$ has been removed. The result is maximally planar and is the third graph represented in Figure \ref{fig:p8b}.
Now suppose edges $(v,w)$ and $(b,y)$ have been removed. The result is still nonplanar, so one more edge must be removed. We obtain a maximally planar subgraph if $(a,x)$, $(v,b)$, or $(v,y)$ is removed. These are the first and second cases in Figure \ref{fig:p8b}. To avoid the first three graphs in Figure \ref{fig:p8b}, we assume that all edges incident to $z$ and all edges incident to $c$ are present, as well as edges $(a,x), (v,b),$ and $(v,y)$.  To obtain a new, maximally planar subgraph, we must remove either $(a,y)$ or $(b,x)$. If we remove both, we obtain a graph isomorphic to a subgraph of the first graph in Figure \ref{fig:p8b}. If we remove just, say, $(a,y)$, then we also must remove edge $(v,a)$ to obtain a planar graph. This result in a subgraph equivalent to a subgraph of the second graph in Figure \ref{fig:p8b}. 

Now we suppose that only one edge of the $Y$ has been removed, and that edge is not incident to $v$. Without loss of generality, we may assume that edge $(c,w)$ has been removed. We note that if the edge $(b,y)$ (or equivalent) is removed, we obtain the fourth maximally planar subgraph pictured in Figure \ref{fig:p8b}. If we instead removed edge $(v,y)$ (equivalently $(v,x)$), we obtain the fifth maximally planar subgraph pictured in Figure \ref{fig:p8b}. If we remove $(c,w)$ and $(b,z)$ (or $(a,z)$), we obtain the sixth graph in Figure \ref{fig:p8b}. Finally, we note that if we removed the edges $(c,w), (c,y), (c,x), (v,b)$ and $(v,a)$, the result would still be non-planar, as we would have $K_{3,3}$ subdivision with the partitions: $\{v,a,b\}$ and $\{x,y,z\}$ (the edge between $v$ and $z$ is subdivided.)  As there are no other possible edges to remove, we have classified all maximal planar subgraphs of $G_8$ that have exactly one $Y$ edge removed.

Now suppose exactly two $Y$ edges have been removed. Suppose first that $(c,w)$ and $(w,z)$ have been removed. Further removing any edge incident to $v$ results in a subgraph of the fifth graph in Figure \ref{fig:p8b}. We similarly do not obtain a maximally planar subgraph if we remove any other edge. 

Now suppose exactly two $Y$ edges have been removed, where one edge is incident to $v$. Without loss of generality, suppose $(v,w)$ and $(c,w)$ have been removed. The only way we obtain a new maximally planar subgraph here is to remove either edge $(v,b)$ or $(v,a)$. This case is represented as the seventh graph in Figure \ref{fig:p8b}.

Finally, a maximal planar subgraph cannot have all three $Y$ edges removed, as this would result in a disconnected graph.

\vskip -1in
\begin{figure}
\hskip .5in\includegraphics[scale=.44]{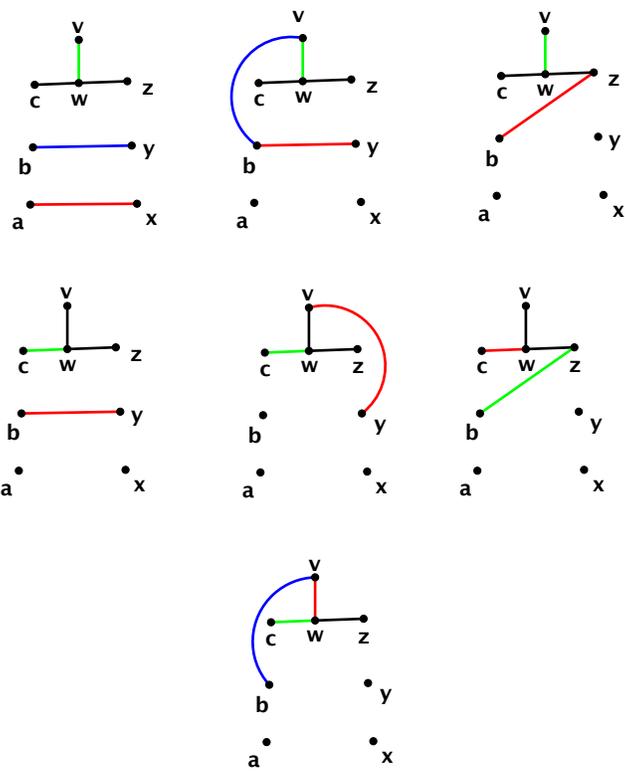}
\vskip -.5in
\caption{The representations of the maximal planar subgraphs of $G_8$ that have removed edges of the $Y$ that would transform $G_8$ back to $K_{3,3,1}$.}
\label{fig:p8b}
\end{figure}

\vskip 1in
Figures \ref{fig:p8second} and \ref{fig:p8bsecond} shows planar embeddings of all these maximal planar subgraphs, together with their unbalanced conflict graphs.

\begin{figure}
\hskip -.35in\includegraphics[scale=.65]{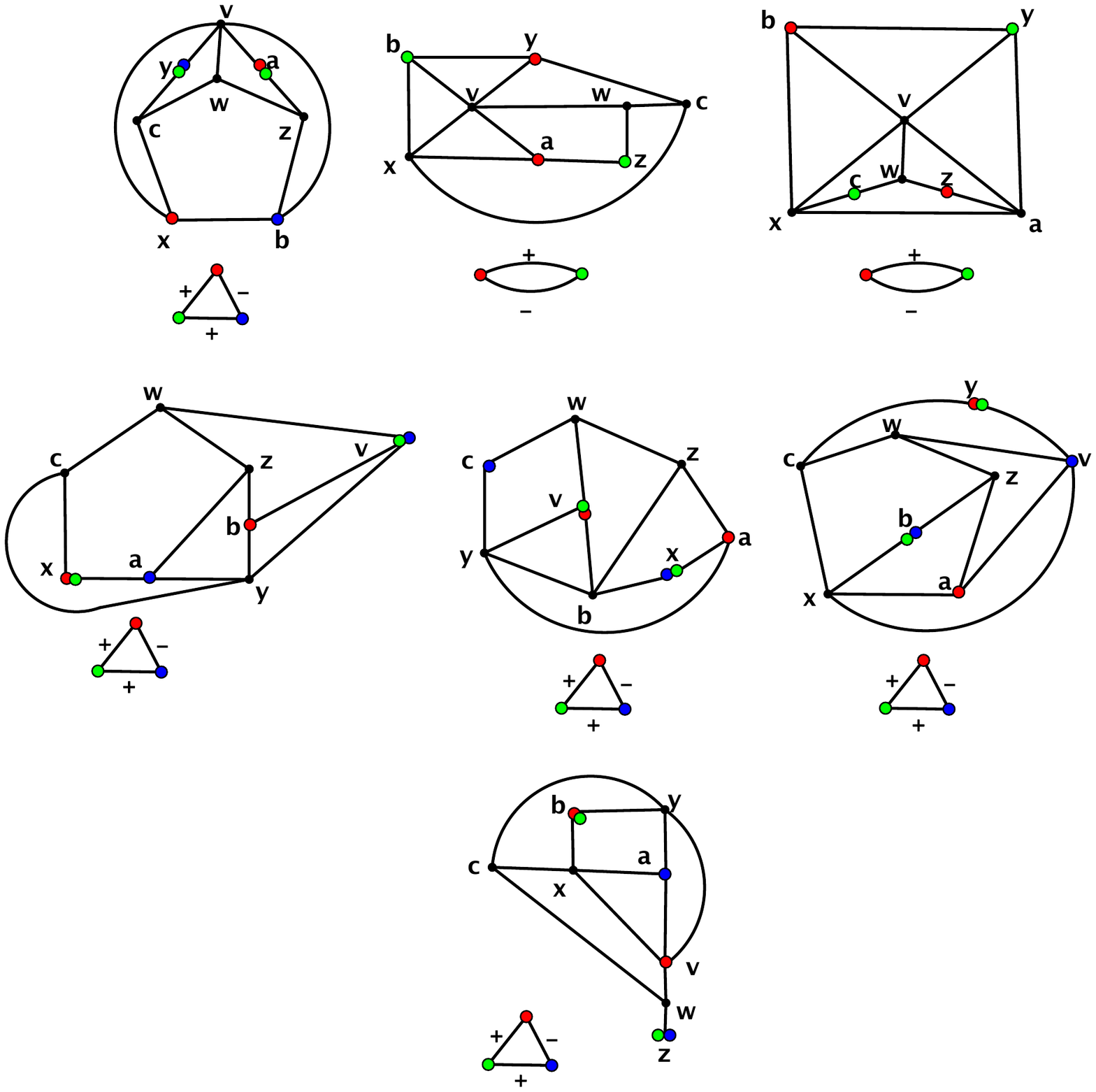}
\vskip -1.15in
\caption{The embeddings and strong conflict graphs of the seven maximal planar subgraphs of $G_8$ from Figure \ref{fig:p8}. Note that the bottom graph pictured is not 2-connected, but the conflict graph is the same regardless of planar embedding.}
\label{fig:p8second}
\end{figure}

\begin{figure}
\hskip -.35in\includegraphics[scale=.65]{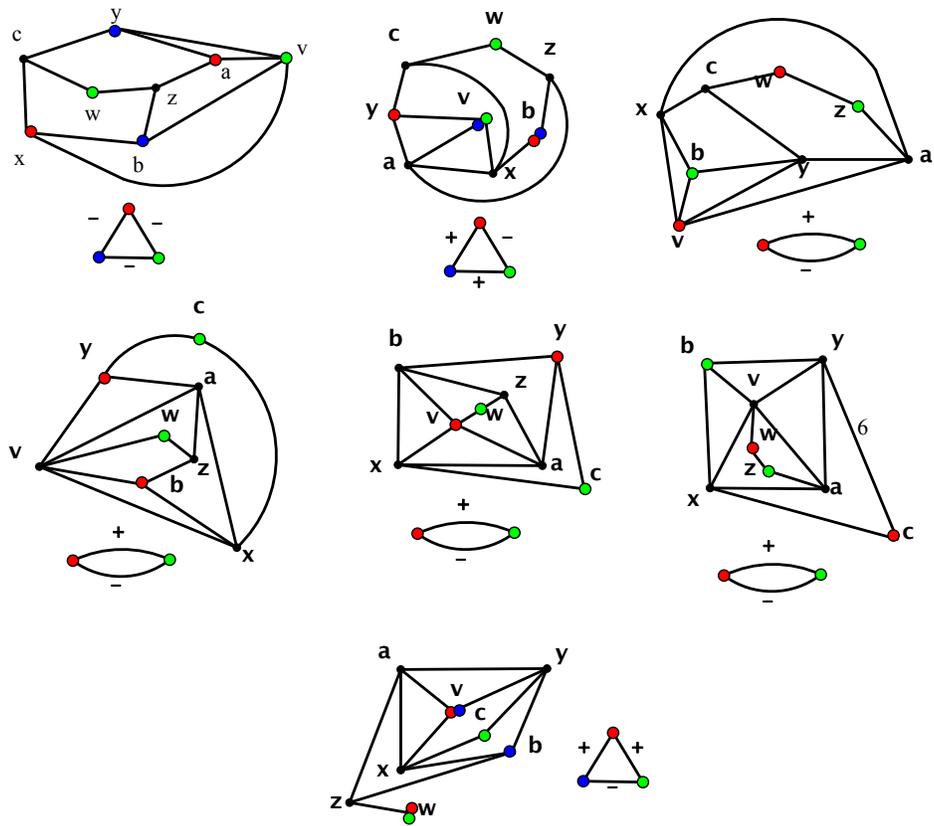}
\vskip -1.4in
\caption{The embeddings and strong conflict graphs of the seven maximal planar subgraphs of $G_8$ from Figure \ref{fig:p8b}. Note that the bottom graph pictured is not 2-connected, but the conflict graph is the same regardless of planar embedding.}
\label{fig:p8bsecond}
\end{figure}
\end{proof}

\begin{figure}
\hskip .25in\includegraphics[scale=.48]{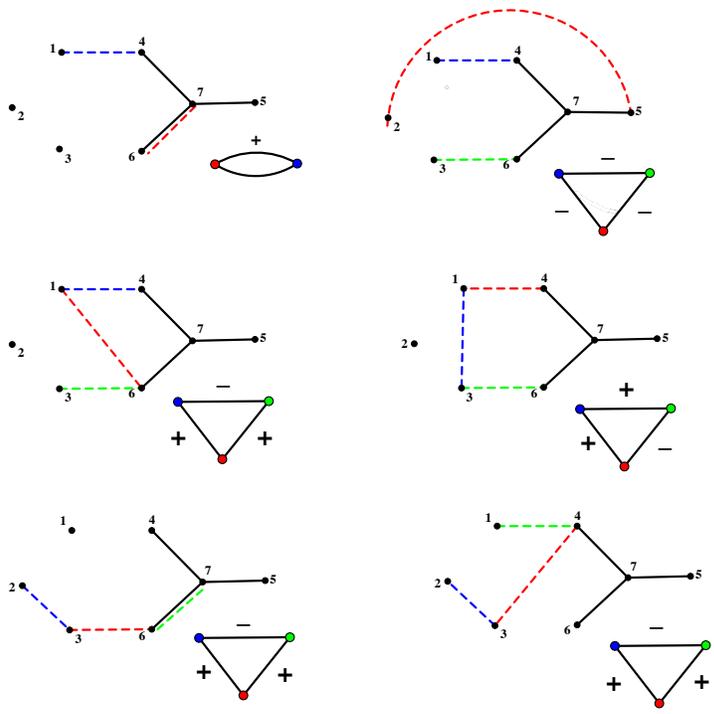}
\vskip -.5in
\caption{Representations of the maximal planar subgraphs of $G_7$.}
\label{fig:p7'}
\end{figure}

\begin{figure}
\hskip .25in\includegraphics[scale=.5]{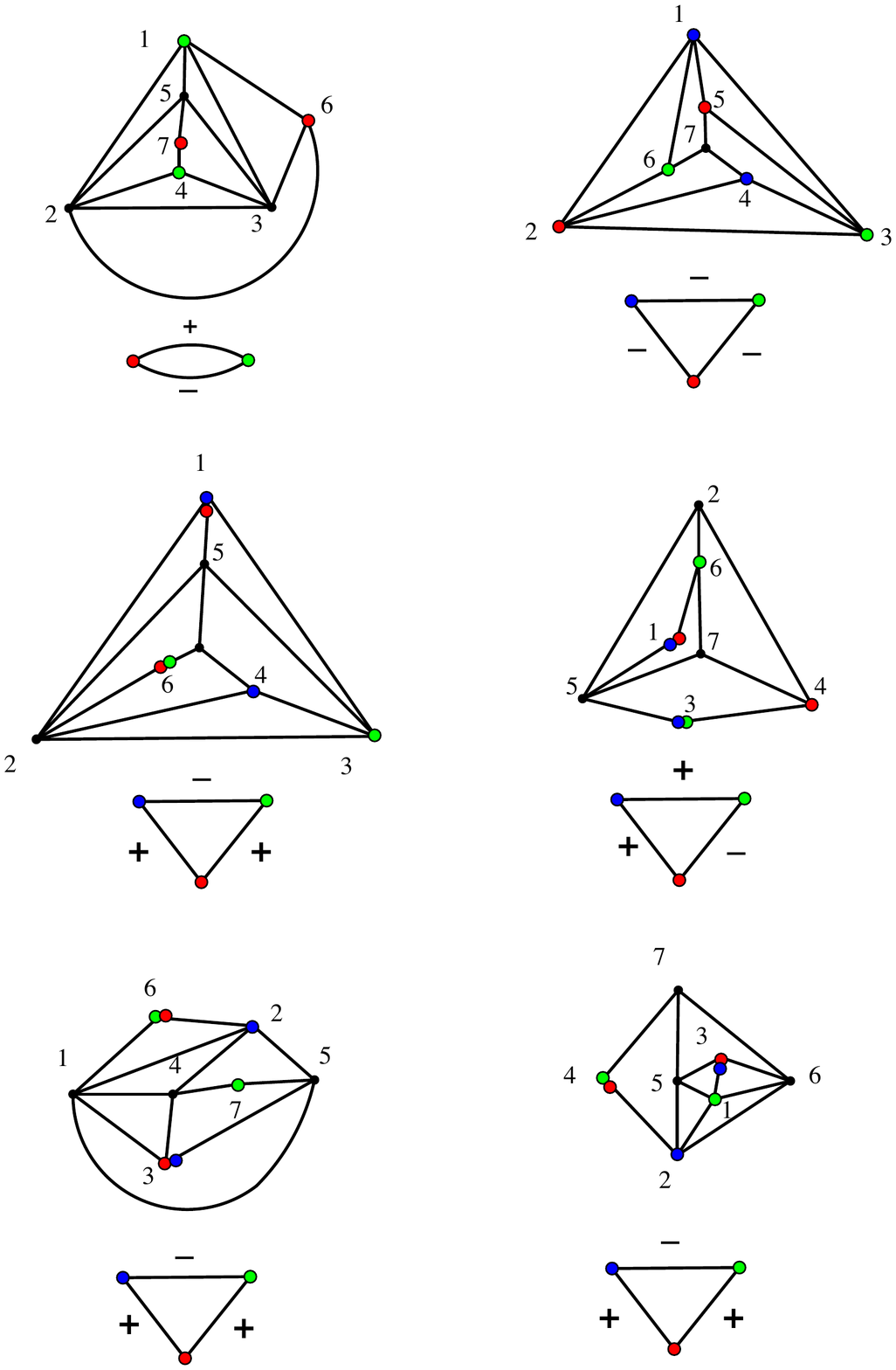}
\vskip -.25in
\caption{The maximal planar subgraphs of $G_7$, with their conflict graphs.}
\label{fig:p7}
\end{figure}

\begin{prop}
The six maximal planar subgraphs of $G_7$ represented in Figure \ref{fig:p7} represent all possible maximal planar subgraphs of $G_7$. Each has an unbalanced strong conflict graph.
\end{prop}
\begin{proof} First, we note that as $G_7$ is obtained from $K_6$ by a single $\Delta-Y$ exchange, we may label the vertices of $G_7$ as $\{1,2,3,4,5,6,7\}$, where vertex $7$ is of degree 3, and it is incident to vertices $4,5$ and $6$. We consider various ways of removing edges to obtain a planar subgraph. First, we consider if edges from the $Y$ have been removed. Note that removing all three $Y$ edges will not result in a maximal planar subgraph, as the result is disconnected. Suppose exactly one edge of the $Y$ has been removed, say edge $(6,7)$. The result is a subdivision of $K_6$ with two adjacent edges removed. It follows that removing an edge incident to $4$ or $5$ results in a maximal planar subgraph (the first one pictured in Figure \ref{fig:p7} and represented schematically in Figure \ref{fig:p7'}). So, we may assume that no edge incident to $4$ or $5$ is removed. We also note that we must remove an edge between $\{1,2,3\}$ and $\{4,5,6\}$, else $K_{3,3}$ will be present as a subgraph. If an edge incident to $6$ is removed, say $(3,6)$, then we must also remove an edge in the cycle $(1,2,3)$. If we remove edge $(2,3)$ or $(1,3)$, we have obtained the fifth graph pictured in Figure \ref{fig:p7}. Suppose we remove edge $(1,2)$. In order to obtain a planar subgraph, we must also remove either $(1,6)$ or $(2,6)$. The result is a subgraph of the fifth graph in Figure \ref{fig:p7}.

If exactly two edges of the $Y$ have been removed, then we still must remove an edge between $\{1,2,3\}$ and $\{4,5,6\}$ to obtain a maximally planar subgraph. In any way we do this, we obtain a subgraph of the first graph in Figure \ref{fig:p7}.

We may now assume that no $Y$ edge has been removed. Given a maximally planar subgraph $G'$, if none of the edges in triangle $(1,2,3)$ have been removed from $G_7$ to obtain $G'$, then the edges removed must form a subgraph of $G_7$ with either 1 or 3 components. Two components are not possible, as $K_6$ with two non-adjacent triangles removed is non-planar, and the $\Delta-Y$ move preserves non-planarity.

To obtain a planar graph, we must remove an edge between $\{1,2,3\}$ and $\{4,5,6\}$. Without loss of generality, suppose that edge $(1,4)$ has been removed. Moreover, we must remove at least two such edges, not both incident to $1$, else there is a $K_{3,3}$ subgraph present, formed by the vertices $\{2,3,7\}$ and $\{4,5,6\}$. By symmetry, we may consider two cases:

Case 1: Edges $(1,4)$ and $(3,6)$ have been removed. If we remove either edge $(2,6)$, or edge $(1,6)$, or edge $(3,4)$, we get either of the second or third maximal planar subgraphs pictured in Figure \ref{fig:p7}. In order for for the removed edges to form a connected or 3 component subgraph, one of these edges must be removed. This concludes Case 1.

Case 2: Edges $(1,4)$ and $(3,4)$ have been removed. If we remove an edge from $\{1,3\}$ to $\{5,6\}$, then, up to equivalence, we get the third planar subgraph in Figure \ref{fig:p7}, so we cannot remove those edges. If edge $(2,4)$ is removed, then we still must remove another edge from $2$ to $\{5,6\}$ to get a planar subgraph. In either case, we get the third planar subgraph in Figure \ref{fig:p7}. There is no other way to remove edges without making the subgraph on the removed edges have exactly 2 components.

We may now assume our planar subgraph has at least one edge from $(1,2,3)$ removed (we still assume no edges of the $Y$ have been removed). We still must remove a graph that has either 3 components or 1 component. The 3 component graph is not possible, given that one of the three edges must come from the triangle $(1,2,3)$. We note that at least 3 edges must be removed, or else by performing a $Y-\Delta$ move, we get a planar subgraph of $K_6$ with fewer than 3 edges removed, which is not possible, as shown earlier. We also note that the removed edges cannot be contained in a triangle, else the resulting graph is nonplanar. In order for $K_{3,3}$ not to be present, at least two edges between $\{1,2,3\}$ and $\{4,5,6\}$ must be removed, not both incident to the same vertex in $\{1,2,3\}$.Without loss of generality, we may assume that edges $(1,4)$ and $(1,3)$ have been removed. We must remove another edge, not incident to $1$. If $(3,6)$ and $(3,5)$ are removed, we get the fourth graph pictured in Figure \ref{fig:p7}. The removal of $(3,4)$ won't result in a planar graph.  If $(2,4)$ is removed, we get the equivalent of the last case in Figure \ref{fig:p7}. The removal of $(2,5)$ or $(2,6)$ won't result in a planar graph. It follows that the only possible maximal planar graphs, for which an edge of $(1,2,3)$ has been removed, are the fourth or sixth cases pictured in Figure \ref{fig:p7}.

Finally, we note that the conflict graph for all of the maximal planar subgraphs of $G_7$ have unbalanced conflict graphs.

\end{proof}

\vskip -.75in

\vskip -.75in

\vskip 1.5in 
\begin{prop} Up to symmetry, there are precisely two maximal planar subgraphs of the classic Petersen graph. Each has an unbalanced strong conflict graph.
\end{prop}
\begin{proof}
Recall that there is a well-known one-to-one correspondence between the vertices of the Petersen graph and two element subsets of $\{1,2,3,4,5\}$, and there is an edge between two vertices if and only if the corresponding two element subsets are disjoint. From this, it readily follows that the symmetry group is $S_5$.  From the upper-left of Figure \ref {fig:PGConflict}, we can see that $PG - \{(\{1,2\},\{4,5\}), (\{3,5\},\{1,4\})\}$ is maximally planar. By symmetry, the following are also maximally planar: (via the permutation $(1,2)(4,5)$: $PG - \{(\{1,2\},\{4,5\}), (\{3,4\},\{2,5\})\}$, via the permutation\\
 $(4,5)$: $PG -\{(\{1,2\},\{4,5\}), (\{3,4\},\{1,5\})\}$, and via the permutation $(1,2)$: \\
$PG - \{(\{1,2\},\{4,5\}), (\{3,5\},\{2,4\})\}$.
 
We can also see that $PG - \{(\{1,2\},\{4,5\}), (\{2,5\},\{1,3\})\}$ is maximally planar, and by symmetry the following are too:  via $(4,5):$\\
$ PG - $ $\{(\{1,2\},\{4,5\}), (\{1,3\},\{2,4\})\}$, via $(1,2)(4,5)$:\\
 $PG - \{(\{1,2\},\{4,5\}), (\{1,4\},\{2,3\})\}$, and via $(1,2)$:\\
  $PG - \{(\{1,2\},\{4,5\}), (\{2,3\},\{1,5\})\}$.

Note that the symmetry given by the permutation $(1,5)(2,4)$ takes $PG - \{(\{1,2\},\{4,5\}), (\{3,5\},\{1,4\})\}$ to $PG - \{(\{1,2\},\{4,5\}), (\{2,5\},\{1,3\})\}$. Note also that given an edge $e$ of $PG$, there are eight edges that are non-adjacent to $e$, and are connected to $e$ by a single edge. There are 2 edges that are not adjacent to $e$ and they are connected to $e$ by a shortest path of length $2$. Further, given two adjacent edges, $e$ and $e'$, the edges that are distance two from $e$ are disjoint from the set of edges that are distance 2 from $e$. For the edge $e=(\{1,2\}, \{4,5\}$, in the previous two paragraphs we have accounted for removing another non-adjacent edge that is connected by a single edge to the original--the result is always a maximal planar subgraph. The other two non-adjacent edges to $(\{1,2\}, \{4,5\})$ are $(\{2,5\},\{1,4\})$ and $(\{1,5\},\{2,5\})$. Removing $e$ and either of these two edges results in a maximally planar subgraph. The permutation $(4,5)$ interchanges either of these two edges paired with $e$. We have a total of two maximal planar subgraphs of PG thus far.

Are there other maximal planar subgraphs of $PG$? By symmetry, we may assume that the first edge removed is $e$ from above. We must remove an edge adjacent to $e$, say $e'$, to get a possibly different maximally planar subgraph. Removing the one edge that is incident to both $e$ and $e'$ will not result in a planar subgraph, as then $PG$ would be apex (that is, becomes planar after the removal of a vertex), and hence not intrinsically linked, which is a contradiction (see \cite{sachs}). It follows that we must remove an edge that is not adjacent to either $e$ or $e'$. By the observation in the previous paragraph, this edge cannot be distance 2 from both. It follows that our resulting graph is not maximally planar.

%

Finally, we see in Figure \ref{fig:PGConflict} that the conflict graph for each of these three maximally planar subgraphs is not balanced.

\vskip -1.5in
\begin{figure}
\hskip .1in\includegraphics[scale=.42]{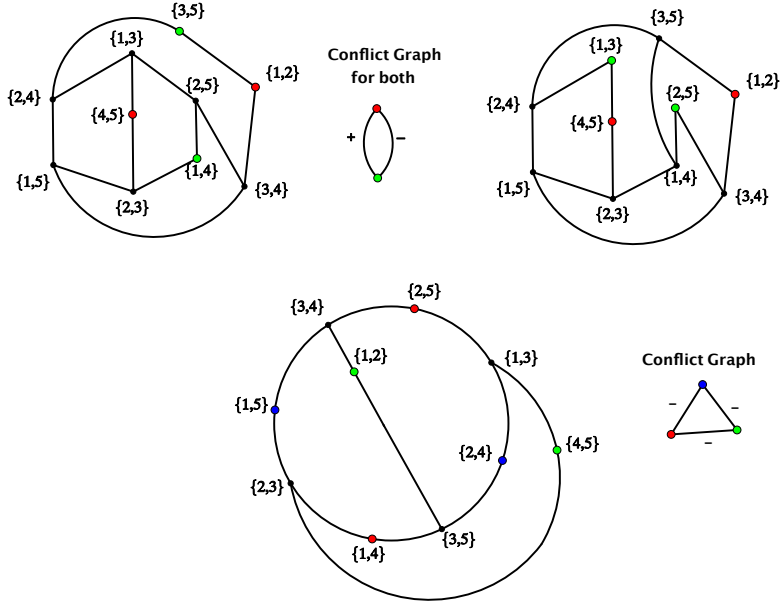}
\vskip -.05in
\caption{The maximally planar subgraphs of the Petersen graph. Note that the top two cases are equivalent. In each case, the associated strong conflict graph is not balanced.}
\label{fig:PGConflict}
\end{figure}

\vskip 1.5in

\end{proof}

\begin{prop} The seven subgraphs represented in Figure \ref{fig:k44} are all of the maximal planar subgraphs of $K_{4,4}-e$. Each has an unbalanced conflict graph.
\end{prop}

\begin{proof}
We denote the vertices of $K_{4,4}-e$ by the two partitions $\{1,2,3,8\}$ and $\{5,6,7,4\}$, where the edge removed is $(4,8)$.

First, we consider maximal planar subgraphs of $K_{4,4}-e$ where no edges of either $Y$ is removed. Such maximal planar subgraphs must result from first removing either three non-adjacent edges or a path of length 3, as $Y-\Delta$ exchanges preserve planarity, and otherwise the result would be a nonplanar subgraph of $K_6$. By checking directly, the two resulting subgraphs are planar, and again they must be maximally planar by the classification of maximally planar subgraphs of $K_6$. This covers the fourth and sixth cases in Figure \ref{fig:k44}.

No suppose that no edges of one $Y$ have been removed, but at least one edge of the other $Y$ has been removed. Any such subgraph, after a $Y-\Delta$ exchange, results in a planar subgraph of $G_7$. By considering the maximal planar subgraphs of $G_7$ that do not have any edges of the triangle that is vertex-disjoint from the $Y$, we may conclude, without loss of generality, that for any such maximally planar subgraph of $K_{4,4}-e$, the edges $(1,4)$ and $(2,5)$ have been removed. To get a planar subgraph, we must remove at least one other edge between $\{1,2,3\}$ and $\{5,6,7\}$. No matter what edge we get, up to symmetry, we get either the first or fifth maximally planar subgraph depicted in Figure \ref{fig:k44}. 

Finally, we suppose that at least one edge of each $Y$ has been removed. To obtain a planar subgraph, we must also remove at least one edge between $\{1,3\}$ and $\{5,6,7\}$. It follows that the resulting graph is equivalent to either the second, third, or last maximally planar subgraph in Figure \ref{fig:k44}.

It follows that, up to symmetry, the maximally planar subgraphs of $K_{4,4}-e$ represented in Figure \ref{fig:k44} are all of the maximally planar subgraphs of $K_{4,4}-e$.

Finally, as Figure \ref{fig:k44conflict} indicates, every maximal planar subgraph of $K_{4,4}-e$ has an unbalanced conflict graph.

\begin{figure}
\hskip .02in\includegraphics[scale=.55]{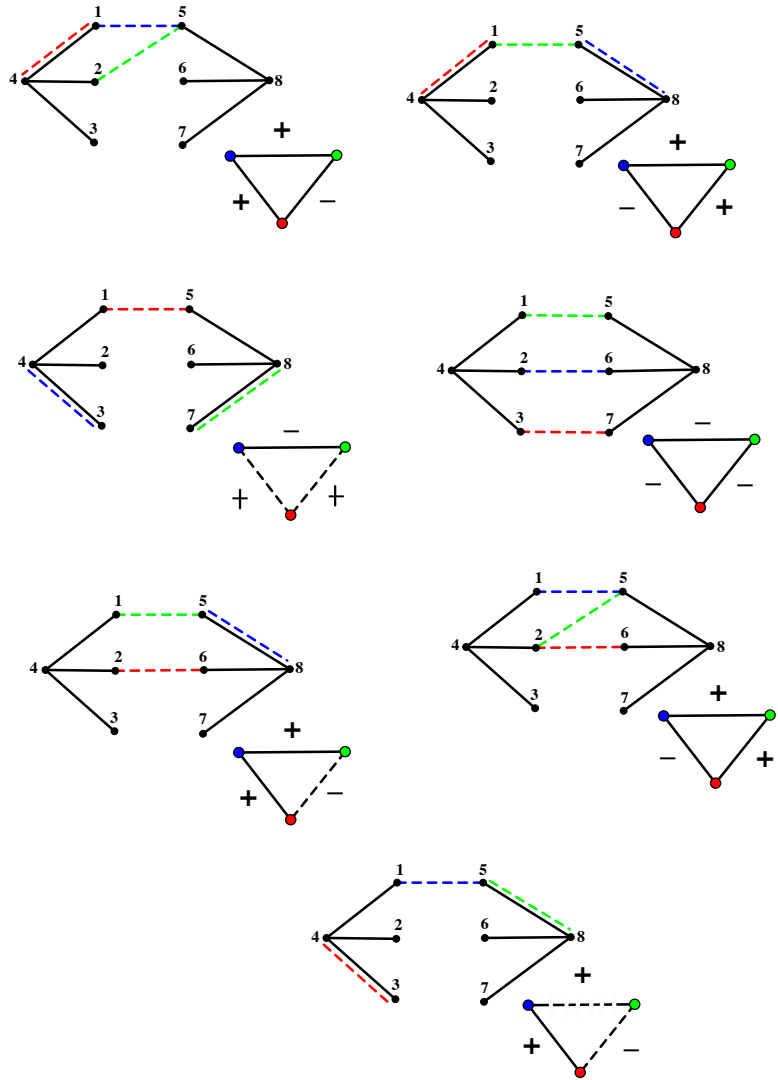}
\vskip -.02in
\caption{Shorthand representations of maximal planar subgraphs of $K_{4,4}-e$.}
\label{fig:k44}
\end{figure}

\begin{figure}
\hskip -.375in\includegraphics[scale=.65]{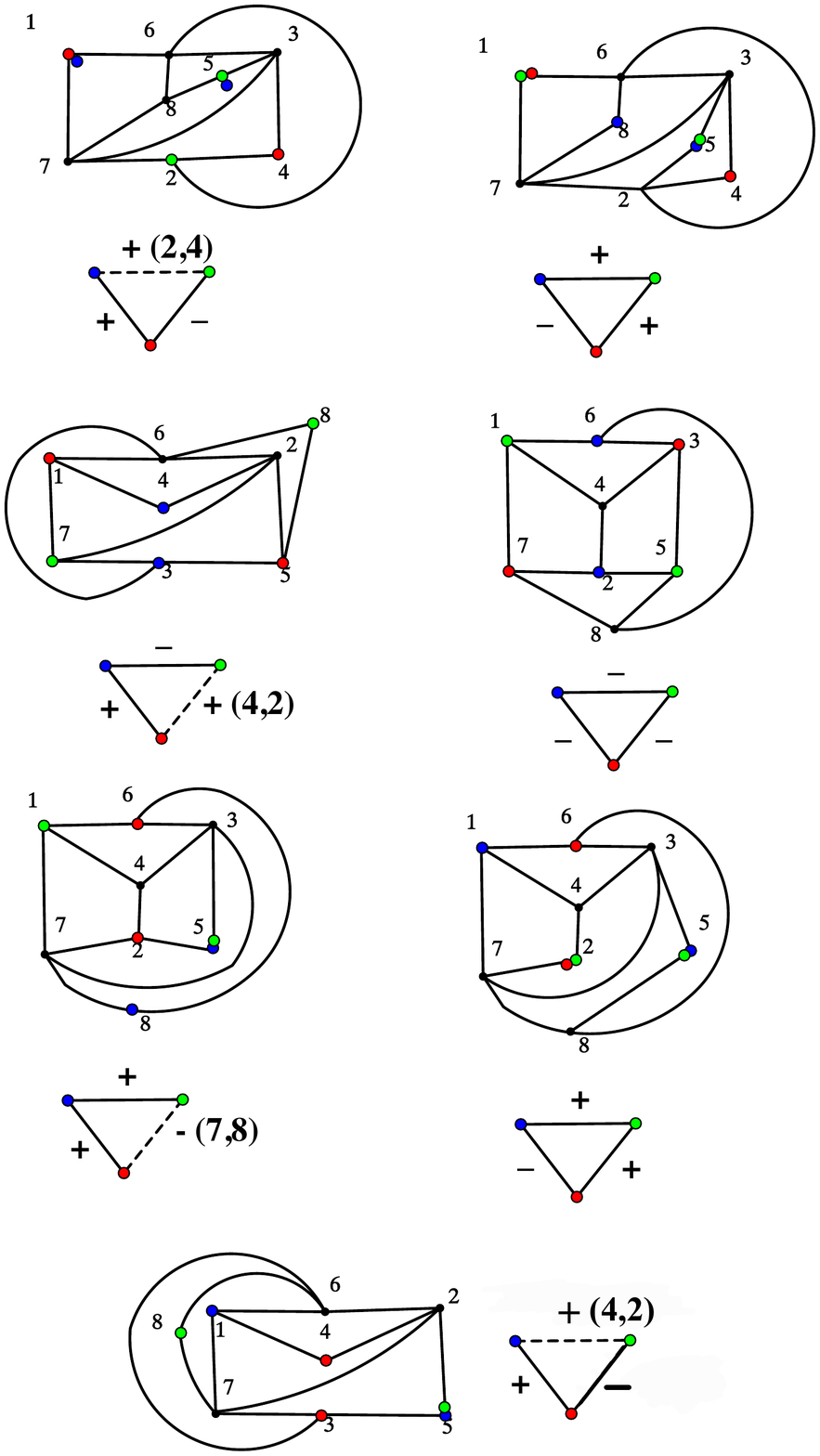}
\vskip -.2in
\caption{Conflict graphs associated to maximal planar subgraphs of $K_{4,4}-e$.}
\label{fig:k44conflict}
\end{figure}

We now justify that the conflict graphs in Figure \ref{fig:k44conflict} are legitimate, in particular the implicit (anti-)conflicts. First, for the top graph, contracting edge $(2,4)$ allows the red fragment to ambient isotope to the face bounded by cycle $(1,6,4-2,7)$, as the red shares endpoints with the other fragments. This allows blue and green to implicitly anti-conflict. 


For the last graph in Figure \ref{fig:k44conflict}, we note first that red and blue strongly anti-conflict, and red and green strongly conflict. As part of the conflict graph process, we build a potentially flat embedding (see \cite{foisy}) starting with the given maximal planar subgraph embedding, and place the red and blue fragments on the same side of the sphere containing the maximal planar subgraph. The green fragment will be placed on the other side of the sphere. We first note that the red and blue fragments can be thought of as forming a tangle from knot theory (see \cite{foisy} proof of Theorem 3.2 for terminology and more references). Since the tangle connects NE to SW and NW to SE, it follows that the tangle formed by the red and blue fragments is ambient isotopic to either the rational tangle +1 or -1 (see \cite{foisy} proof of Theorem 3.2 and Figure \ref{tangle}). 



\begin{figure}

\hskip 1in\includegraphics[scale=.75]{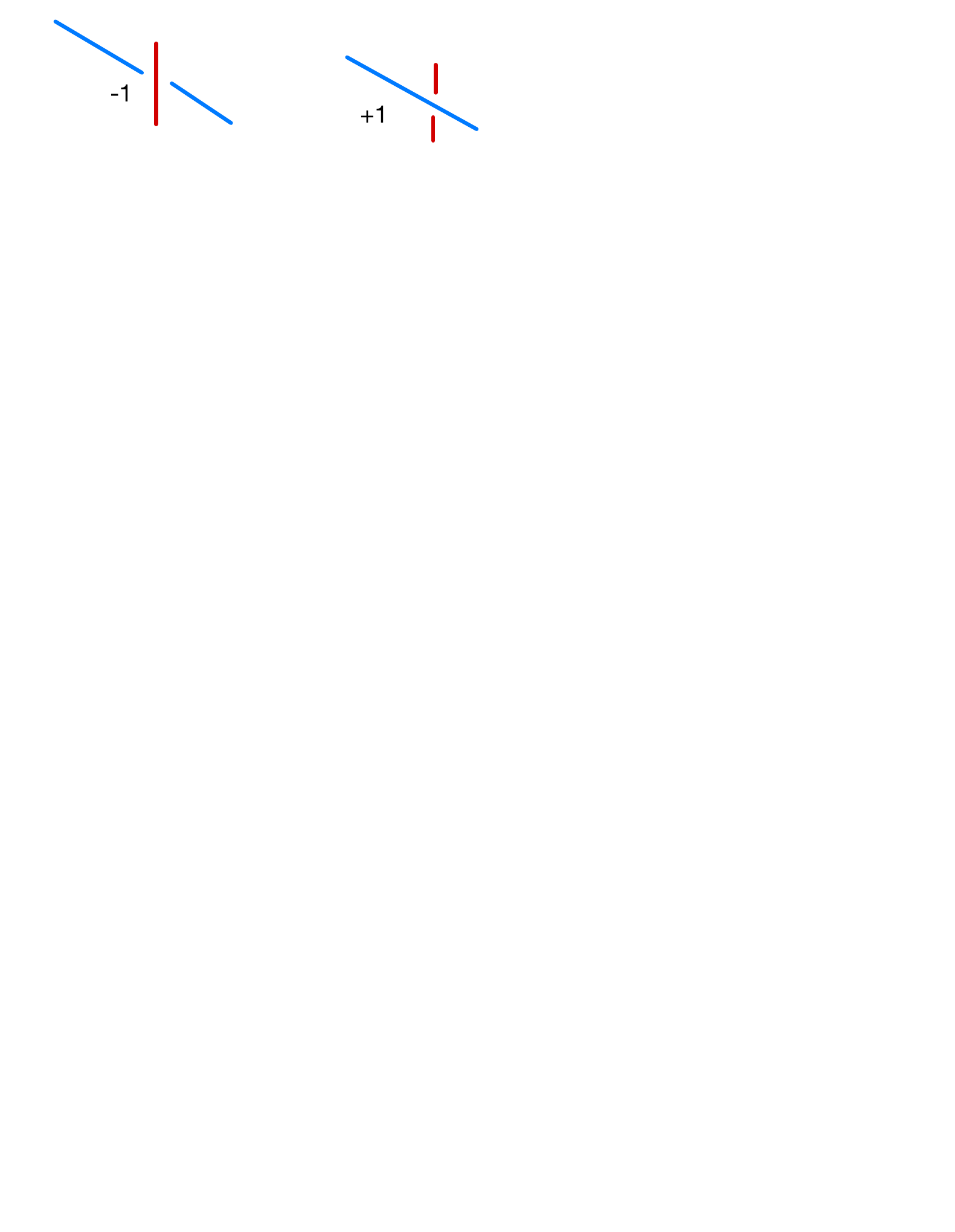}
\vskip -5.3 in
\caption{The two tangles possible for a potentially flat embedding.}
\label{tangle}
\end{figure}

If blue and red form +1, then by contracting edge $(4,2)$, the red fragment can slide onto the sphere (into the face bounded by cycle $(2, 7, 3, 5)$) so that blue and green anti-conflict. If blue and red form -1, then by contracting edge $(4, 2)$, red slides into the face bounded by $(2, 7, 3, 5)$ and again blue and green anti-conflict. Since blue and green anti-conflict for every possible potentially flat embedding with red and blue on the same side and blue and green on the opposite side, then blue and green implicitly anti-conflict.

For the third graph, for a potentially flat embedding, red and blue are on the same side, with green on the other side. Similar to the last graph, red and blue form a tangle that must be equivalent to either +1 or -1. If the tangle is +1, then contracting edge $(4,2)$ allows blue to land in the face bounded by $(7, 2, 5, 3)$ and red and green anti-conflict. If the tangle is -1, then a nonsplit link is present: $(4,2,3,7)$ and $(1, 5, 8, 6)$. As we have considered all possible potentially flat embeddings, this shows that red and green implicitly anti-conflict.

For the fifth graph, in a potentially flat embedding, the red and blue will be on the same side (as well as green). Similar to the last case and third case, the red and blue form a tangle that can be closed up two ways: connecting via paths $(5, 3,6)$ and $(2,7,8)$ or through $(2,5)$ and $(8,6)$. Again, in a potentially flat embedding, the tangle must be equivalent to +1 or -1 (for consistency, we will consider the green fragment as having endpoints at NW and SE). Then +1 will be the tangle for which, after contracting edge $(7,8)$, blue can deform to the face bounded by $(2,5,3,7)$. This contraction yields a conflict between green and red. The link being equivalent to $-1$ would have guaranteed a nonsplit link present involving blue, red, and the maximal planar graph: namely $(1,6,2,4)$ and $(5,3,7,8)$. It follows that green and red implicitly anti-conflict.
\end{proof}
\begin{prop}
The six graphs in Figure \ref{fig:P9first} represent, up to isomorphism, all possible maximal planar subgraphs of $G_9$ with exactly two edges removed.
\end{prop}

\begin{proof}
We label the vertices of $G_9$ as $\{0,1,2,3,4,5,6,7,8\}$, where $(0,7,8)$ forms the only triangle, and if that triangle is removed, and the resulting degree 2 vertices are smoothed out, the resulting $K_{3,3}$ has $\{1,3,5\}$ in one partition and $\{2,4,6\}$ in the other. We first consider maximal planar subgraphs that have exactly two edges removed. Note that we cannot remove exactly two triangle edges and obtain a planar graph. If one triangle edge is removed, say $(0,7)$, then the only way to obtain a planar graph is to remove either $(1,4)$ or $(2,3)$ (see the fourth graph in Figure \ref{fig:P9first}). If we remove an edge adjacent to $(0,7)$, then the result has to be nonplanar. The only remaining options are to remove, say $(3,6)$ (contract $(3,7)$ and $(8,6)$ to get $K_{3,3}$) or to remove $(6,8)$, in which case we can contract edge $(0,1)$ and $(3,7)$ to get $K_{3,3}$ with partitions $\{1,3,5\}$ and $\{2,4,8\}$.

Now suppose no triangle edge has been removed. Suppose edge $(0,1)$ has been removed first. If we remove $(2,3)$ (equivalently $(2,5)$), $(3,6)$ (equivalently $(4,5)$) or $(4,7)$ (equivalently $(6,8)$), we obtain a maximally planar subgraph (see the first three graphs in Figure \ref{fig:P9first}).  Let us consider what other non-adjacent, non-triangle, edges can be removed. If we remove $(3,7)$ (equivalently $(5,8)$), we obtain $K_{3,3}$ by contracting $(0,7)$. The two partitions are $\{5,6, (0,7)\}$ and $\{2,4,8\}$ (note that ``edges" from $6$ pass through the vertex $1$).
There are no other possible non-adjacent, non-triangle edges that we can remove.

Finally, suppose no triangle edges are removed, nor any edges incident to the triangle. Suppose, without loss of generality, that edge $(1,4)$ has been removed. If we remove edge $(2,3)$, the result is planar. This is also true if we remove the edge $(2,5)$ (equivalently, $(3,6)$)--see the last two graphs shown in Figure \ref{fig:P9first}. Thus, we have shown every maximally planar subgraph of $G_9$ with exactly two edges removed. We note that the strong conflict graph for each is an unbalanced bigon.\\

\vskip -.25in
\begin{figure}
\hskip .1in\includegraphics[scale=.45]{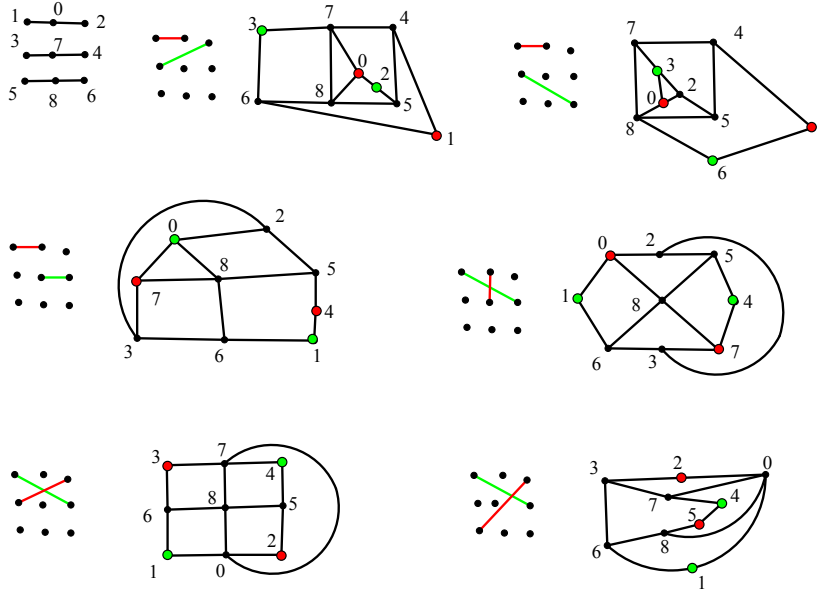}
\vskip -.2in
\caption{The maximal planar subgraphs of $G_9$ with exactly two edges removed.}
\label{fig:P9first}
\end{figure}

\end{proof}

\begin{prop}
The four graphs in Figure \ref{fig:P9second} represent, up to isomorphism, all possible maximal planar subgraphs of $G_9$ with exactly three edges removed. Moreover, the 10 graphs pictured in Figures \ref{fig:P9first} and \ref{fig:P9second} represent all possible maximal planar subgraphs of $G_9$.
\end{prop}

\begin{proof}
We first note that Figure \ref{fig:P9second} shows four maximally planar subgraphs of $G_9$ with three edges removed. 

Suppose a maximal planar subgraph of $G_9$ resulted from removing all edges of the triangle. In order to obtain a planar subgraph, an additional edge must be removed. If a 3-3 (where 3 refers to the degree of the vertices in $G_9$) edge is removed, then we obtain a subgraph of the fourth graph pictured in Figure \ref{fig:P9first}. If a 3-4 edge is removed, we again obtain a subgraph of the third graph pictured in Figure \ref{fig:P9second}. Thus, a maximal planar subgraph must contain at least one triangle edge. 

Suppose exactly two triangle edges have been removed, say $(0,7)$ and $(7,8)$. If we remove a $3-3$ edge, then again, we obtain the fourth graph in Figure \ref{fig:P9first}, unless the edge removed is $(1,6)$ or $(2,5)$. If both are removed, we obtain a subgraph of the fifth graph in Figure \ref{fig:P9first}. Otherwise, if we remove only one such edge, the result in non-planar. If we remove a 3-4 edge, we obtain the third graph in Figure \ref{fig:P9second}, unless we remove $(3,7)$ or $(7,4)$. If we do remove $(0,7), (7,8), (3,7), $ and $(7,4)$, the result is non-planar. If we remove both a $3-3$ and a $3-4$ edge, say, both $(4,7)$ (or $(3,7)$) and $(2,5)$ (or $(1,6)$), then we obtain the equivalent of a subgraph of the second graph in Figure \ref{fig:P9first}.

By the argument in the previous paragraph, we may assume that no more than 1 triangular edge has been removed.
\vskip -.15in
\begin{figure}
\hskip .2in\includegraphics[scale=.4]{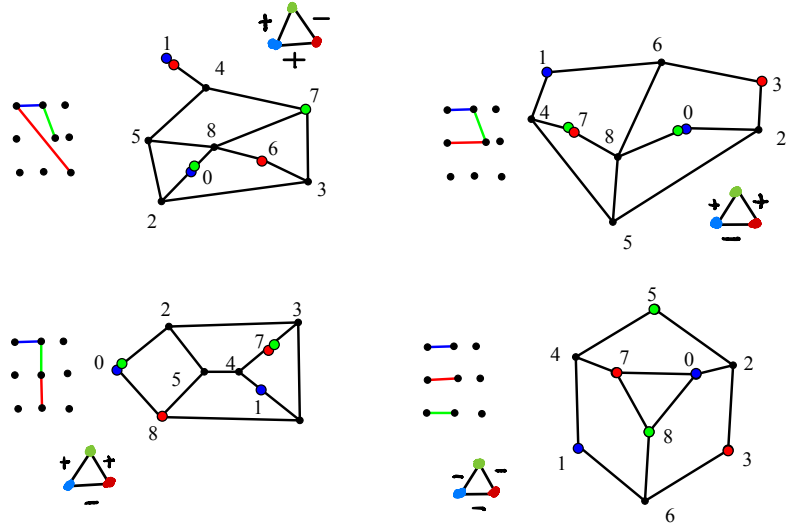}
\vskip -.35in
\caption{The maximal planar subgraphs of $G_9$ with exactly three edges removed. Note the conflict graph for the first would be the same even if re-embedded. }
\label{fig:P9second}
\end{figure}

\vskip .15in

Now we suppose that exactly one triangle edge has been removed, say $(0,7)$. If we further removed $(1,4)$ or $(2,3)$, we would obtain the fourth graph pictured in Figure \ref{fig:P9first}. Suppose no edge adjacent to $(0,7)$ has been removed. If we only remove edge $(3,6)$ (equivalently, $(4,5), (1,6)$ or $(2,5)$), then the result is non-planar. If we remove another edge non-adjacent to $(0,7)$, we obtain a subgraph of either the fourth, fifth or sixth graph in Figure \ref{fig:P9first}. So now we assume that $(0,7)$ and an additional edge has been removed, with the removed edge adjacent to $(0,7)$ (and not from the triangle). Without loss of generality, suppose edge $(1,0)$ has been removed. The result of removing $(0,7)$ and $(1,0)$ is nonplanar, so at least one more edge must be removed. If edge $(3,7)$ is removed, we obtain the second graph pictured in Figure \ref{fig:P9second}. Removing edge $(0,2)$ results in a disconnected graph. Removing $(7,4)$ or $(8,6)$ results in a subgraph of the third graph in Figure \ref{fig:P9first}. Removing any edge between two degree 3 vertices results in a subgraph of either the first, second, or third graphs in Figure \ref{fig:P9first}, or the first in Figure \ref{fig:P9second}. Finally, removing only edge $(5,8)$ results in a graph that is still nonplanar with partitions $\{2,6,7\}$ and $\{3, (8,0), (1,4,5)\}$ forming a $K_{3,3}$ subdivision.

Finally, we suppose that no triangle edges have been removed. To get a new maximally planar subgraph, at least 3 edges must be removed, and they cannot all be incident to the same vertex (else the result is nonplanar). If we remove any two vertex disjoint 3-3 edges, the result is either the fifth or sixth graph in Figure \ref{fig:P9first}. Similarly, if we remove disjoint $3-4$ and $3-3$ edges, the result is either the first or second graphs pictured in Figure \ref{fig:P9first}. If we remove exactly two disjoint $3-4$ edges, either we get the third graph in Figure \ref{fig:P9first}, or we get a nonplanar graph. Without loss of generality, suppose $(0,1)$ and $(3,7)$ have been removed. The resulting graph is non-planar. If we remove another $3-4$ edge, we get either a subgraph of the third graph in Figure \ref{fig:P9first}, or a subgraph of the last graph in Figure \ref{fig:P9second}.  If we remove another $3-3$ edge, we get either a subgraph of the first or second graphs in Figure \ref{fig:P9first}.

As this covers all cases, we have shown that the graphs in Figure \ref{fig:P9first} and Figure \ref{fig:P9second} represent all possible maximally planar subgraphs of $G_9$. We note that the associated strong conflict graphs are all unbalanced.
\end{proof}

\section*{Acknowledgments}
This work started as an REU project with students Cara Nickolaus, Justin Raimondi, Joshua Wilson and Liang Zhang. We wish to thank the SUNY Potsdam-Clarkson University REU Program and the NSA and NSF for financial support under NSA Grant H98230-11-1-0206 and NSF Grant DMS-1004531. We would like to thank  Dr. Tom Zaslavsky and Dr. Garry Bowlin for helpful conversations. 

 \bibliographystyle{elsarticle-num} 
 \bibliography{cas-refs} 

\end{document}